\documentclass{article}

\usepackage[T1]{fontenc}
\usepackage[utf8]{inputenc}
\usepackage{lmodern} %sorgt für schöne Lettern
\usepackage[english]{babel}
\usepackage{hyperref}
\usepackage{authblk}

\usepackage{amssymb}
\usepackage{amsmath}
\usepackage{mathtools}
\usepackage{bm}
\usepackage{xfrac}
\usepackage{graphicx}
\usepackage{subcaption}
\usepackage{xcolor}
\usepackage{siunitx} % for SI units
\usepackage{cleveref} % for Cref

\usepackage{algorithm}
\usepackage[noend]{algpseudocode}
\algnewcommand{\IfThenElse}[3]{% \IfThenElse{<if>}{<then>}{<else>}
	\State \algorithmicif\ #1\ \algorithmicthen\ #2\ \algorithmicelse\ #3}
\newcommand{\algorithmicbreak}{\textbf{break}}

\usepackage{svg} % for SVG support
\svgpath{{svgs/}}

\usepackage{tikz} % all this block for tikz support
\usetikzlibrary{shapes.geometric,calc,positioning}
\usepackage{tikz-cd}
\usepackage{tikz-3dplot}
\usepackage{pgfplots}				% tikz möglich
\usepackage{varwidth}
\newlength\figureheight 
\newlength\figurewidth
\newlength\svgwidth
\newlength\boxwidth
\pgfplotsset{compat=newest,grid style={dotted,black}}
\usepgfplotslibrary{external} 			% benötigt --shell-escape im Ausgabeprofil
\usepgfplotslibrary{patchplots}
\graphicspath{{svg-inkscape/},{tikzs/}}

\newcommand{\tunderbrace}[2]{\underbrace{#1}_{\textstyle#2}}
\newcommand{\R}{\mathbb{R}}
\newcommand{\dd}{\,\mathrm{d}}
\newcommand{\crack}{\Gamma}
\newcommand{\domain}{\mathcal{B}}
\newcommand{\displacement}{\bm{w}}
\newcommand{\strain}{\bm{\epsilon}}
\newcommand{\stress}{\bm{\sigma}}
\DeclareMathOperator{\Div}{div}
\DeclareMathOperator{\trace}{tr}

\begin{document}

%\begin{frontmatter}



\title{Crack propagation in anisotropic brittle materials: from a phase-field model to a shape optimization approach}


\author[1]{Tim Suchan\footnote{Corresponding author: suchan@hsu-hh.de}}
\author[1]{Chaitanya Kandekar}
\author[1]{Wolfgang E. Weber}
\author[2]{Kathrin Welker}
\affil[1]{Professorship of Structural Analysis, \textsc{Helmut Schmidt} University / University of the Federal Armed Forces Hamburg, Holstenhofweg 85, 22043 Hamburg, Germany}
\affil[2]{Professorship of Mathematical Optimization, Technische Universität Bergakademie Freiberg, Akademiestraße 6, 09599 Freiberg, Germany}
% \affiliation{organization={},addressline={}, city={}, postcode={}, state={},country={}}

\maketitle

\begin{abstract}
%%1. Quick overview of the topic
%===========================================
%Common methods to simulate fracture propagation include for example the insertion of cohesive elements {in the mesh at the {expected} fracture area} or a regularization with a phase-field method. %~\cite{Miehe2010}.
%
The phase-field method is based on the energy minimization principle which is a geometric method for modeling diffusive cracks that are popularly implemented with irreversibility based on \textsc{Griffith}'s criterion. %~$G_c$.
%===========================================
%%2. Clearly mention the research questions
%===========================================
%This method yields a strong mesh dependence of the fracture path due to {the} {user-defined position of} {the} {inserted} {cohesive} {elements}. 
This method requires a length-scale parameter that smooths the sharp discontinuity, which influences the diffuse band and results in mesh-sensitive fracture propagation results.
%===========================================
%3. Mention the research methodology
%===========================================
Recently, a novel approach based on the optimization on \textsc{Riemann}ian shape spaces has been proposed, % in~\cite{Suchan2023}, 
where the crack path is realized by techniques from shape optimization. This approach requires the shape derivative, which is derived in a continuous sense and used for a gradient-based algorithm to minimize the energy of the system. Due to the continuous derivation of the shape derivative, this approach yields mesh-independent results.
In this paper, the novel approach based on shape optimization is presented, followed by an assessment of the predicted crack path in anisotropic brittle material using numerical calculations  from a phase-field model.
\end{abstract}

%===========================================
%% Graphical abstract
% \begin{highlights}
% \end{highlights}
%===========================================
%\begin{keyword}

%\KW{@all: Tim and I have rewritten Abstract, Introduction, subsection 2.1,     . Moreover, we have changed the title.}
%
%\KW{Notes: We need to check the following points before submitting. In view of our tough deadline (Tuesday before 6pm submitting to ArXiv), I will add Names to the following ToDos: \TS{I am too tired to continue now. Please someone do the remaining 2 tasks tomorrow in the morning @CK/WW/KW}
%\begin{itemize}
%\item subscripts and superscripts should be italic or not but consistent \TS{done}
%\item transposed vectors or matrices should be $^\top$ and not $^T$ \TS{done}
%\item We should not use formulations in the simple future (will) in the paper.  \TS{done}
%\item Only numeration of formulas cited in the paper.
%\item "colon" instead of double-dot for definitions (TS) \TS{done}
%\item \TS{Check References for consistent formatting}
%\end{itemize}
%}

\section{Introduction}\label{sec:Introduction}
%=======================================================================
%% General Guideline for Introduction:
%1. establishing the context, background and/or importance of the topic
%2. identifying a problem, controversy or a knowledge gap in the field of study
%3. giving a brief review of the relevant academic literature
%4. stating the aim(s) of the research and the research questions or hypotheses
%5. providing a synopsis of the research design and method(s)
%6. explaining the significance or value of the study
%7. providing an overview of the dissertation or report structure
%=======================================================================

% Literature part for Phase-field method
Phase-field modelling is a widely used approach to solve various types of multi-physics phenomena, particularly in fracture problems, in both the physics and the mechanics community~\cite{ambati15}.
Application include mechanical engineering, civil engineering, hydraulic engineering~\cite{Mikota2003} or aerospace engineering~\cite{Hoeche2021, Rauter2023, Rottmann2023}.
In general, the phase-field method is a technique which is used to model sharp interfaces with the help of a continuous scalar field variable, which facilitates the difference between multiple physical phases. 
In fracture mechanics applications, multiple phases are interpreted as solid (intact) or fully broken. 
Thus, the phase-field approach is used for modelling brittle fracture to predict crack initiation, propagation, merging, and branching.
Due to its simplicity in implementation, this methodology has gained wide interest in the engineering community since 2008. 
%Since then, many scientists have worked in this field and developed phase-field approaches for finite elements, isogeometric analysis, and lately also for the virtual element technology.\TS{suggestion: remove previous sentence because the same is mentioned in the following sentences}
The main driving force for these developments is the possibility to handle complex fracture phenomena within numerical methods in two and three dimensions. 
In recent years, several brittle fracture applications~\cite{alessi2018,kakouris17, kuhn10, kuhn15, miehe+welschinger+hofacker10, singh16, zhang17} and an extension to anisotropy~\cite{bleyer+alessi18,Nguyen2017a, Fadianiso17} have been published. It was further developed also for isogeometric analysis~\cite{borden+hughes+landis+verhoosel14} and the virtual element method (VEM)~\cite{aldakheel+blaz+wriggers18}.
Applications for ductile fracture are also described in~\cite{Aldakheel2021, ambati+gerasimov+lorenzis15}.
Furthermore, the phase-field method is also applied in multi-physics problems such as corrosion~\cite{Aldakheel2022, Mai2016} and corrosion cracking~\cite{Cui2021, Cui2022, KandWebHoech:2022}.
These advances are mainly motivated by the capacity to handle complex fracture propagation in two- and three-dimensional numerical approaches.  
These investigations include a wide range of topics, including the modeling of 2D and 3D small and large strain deformations, variational formulations, multi-scale %/physics 
problems, mathematical analysis, various decompositions, and discretization approaches with numerous applications in science and engineering.

% Literature part for shape optimization
Shape optimization has seen frequent uses for medical \cite{Naegel2015}, structural \cite{Loehner2003,Othmer2008, Soto2004}, fluid-mechanical \cite{Allaire2004,Schulz2016,Upadhyay2021}, and acoustical applications \cite{Kapellos2019, Schmidt2016}
%Furthermore, shape optimization has also been employed in the context of multi-physics problems, cf., e.g., \cite{Bletsos2021}. Gradient-based optimization algorithms have been formulated and are often used \cite{Schulz2016}---however, other shape optimization algorithms exist, e.g., for gradient-free approaches \cite{Baeck1996} or approaches based on second-order derivatives \cite{Hintermueller2004,Schulz2014,Schmidt2018}. 
but the usage of shape optimization for fracture problems is relatively new. Recently, a novel method for solving brittle isotropic fracture problems by using a shape optimization approach that minimizes the energy in the system by a part of gradient-based techniques has been proposed in~\cite{Suchan2023}. In that publication, the fracture is described by a surface of the discretized domain, and therefore has sharp edges, in contrast to a phase-field approach. Other shape optimization algorithms exist, e.\,g., gradient-free approaches \cite{Baeck1996} or approaches based on second-order derivatives \cite{Hintermueller2004,Schmidt2018, Schulz2014}, but show other disadvantages regarding the number of partial differential equation (PDE) evaluations or existence of and analytical effort of calculating higher derivatives. Due to the novel nature of the approach, several challenges are still present also for a gradient-based approach: extension of the approach to different material models, the implementation of the irreversibility of the fracture, and the computation of the deformation of the computational domain from the derivative of the energy with respect to the shape. In this manuscript, we aim to extend the gradient-based approach  for brittle fracture in~\cite{Suchan2023} to anisotropic material behavior and establish a comparison between this approach and the previously-mentioned phase-field method.
%
%This methodology is exploited in the numerical implementation by solving brittle fracture problems and comparing it to the phase-field method.
The fracture response of the isotropic and aniostropic elastic material is analyzed in two boundary value problems to illustrate the shape optimization approach in applications for both isotropic and anisotropic material.

The paper is organized as follows: 
\Cref{sec:Modelling} outlines the brittle fracture modelling and the fundamentals of both phase-field and shape optimization methods, to describe the underlying governing equations along with the necessary assumptions to simplify the problem.
\Cref{sec:Implementation} includes a distinct structure of balance equations and their algorithmic implementation to calculate the crack propagation in anisotropic elastic medium using the finite element method (FEM) with shape optimization and phase-field approaches. 
Some numerical tests are carried out to substantiate our algorithmic developments and to test effects of anisotropy on the crack path in \Cref{sec:NumericalExamples}. Herein, the results of numerical tests are discussed. 
Finally, the summary and conclusions are provided in \Cref{sec:SummaryConclusions}.
	
%=======================================================================
\section{Modelling}
\label{sec:Modelling}
%=======================================================================

In this section, we first outline a commonly-used model for brittle fracture in \Cref{subsec:FractureModel}. Then,  \Cref{subsec:AdaptionFractureModelPhaseFieldSetting} presents how this general model can be adapted to a phase-field setting. Lastly, \Cref{subsec:AdaptionFractureModelShapeOptimizationSetting} describes the usage of the model for brittle fracture in a shape optimization algorithm.

%-----------------------------------------------------------------------
\subsection{Fracture model}
\label{subsec:FractureModel}
%-----------------------------------------------------------------------

As is well-known (cf., e.g., \cite{Francfort1998}), the stored energy due to linear elastic deformation of a domain $\domain\subset \R^2$ omitting any volume and surface loads can be computed as
\begin{align*}
	E_{bulk} = \inf_{\displacement \in H^1(\domain, \R^2)} \int_\domain \Psi_0(\strain(\displacement)) \dd\bm{x},
\end{align*}
where $\Psi_0$ denotes the energy density function and is defined as 
\begin{align*}
	\Psi_0(\strain(\displacement)) = \frac{1}{2} \strain(\displacement) : \mathbb{C} : \strain(\displacement) = \frac{1}{2} \strain^V(\displacement) \cdot \bm{C} \cdot \strain^V(\displacement).
\end{align*}
Herein, the strain tensor is defined as $\strain (\displacement)= \frac{1}{2} \left(\nabla \displacement + {\nabla \displacement}^\top \right)$.
The material stiffness tensor $\mathbb{C}$ and $\bm{C}$ describe the material behavior, and the stress is defined as $\stress(\displacement) = \mathbb{C} : \strain(\displacement)$ and $\stress^V(\displacement) = \bm{C} \cdot \strain^V(\displacement)$, respectively. Here, in view of the implementation of anisotropy of the material stiffness, strain, and stress tensor have been expressed in \textsc{Voigt} notation\footnote{By using the symmetry of the tensors, the \textsc{Voigt} notation \cite{Voigt1910} describes canonical isomorphisms 
	\begin{align*}
	\R_{sym}^{2 \times 2 \times 2 \times 2} & \to \R^{3 \times 3},\, \mathbb{C} \mapsto \bm{C},\\
	\R_{sym}^{2 \times 2}& \to \R^{3},\, \strain \mapsto \strain^V, \text{ and}\\
	\R_{sym}^{2 \times 2} &\to \R^{3},\, \stress \mapsto \stress^V.
	\end{align*}
}. 
In case of isotropic behavior, $\bm{C}=\bm{C}_0$. 
For anisotropic materials, $\bm{C}(\theta) =\bm{P}(\theta) ^\top \bm{C}^{ref} \bm{P}(\theta)$ with a rotation matrix defined for a rotation angle $\theta \in [0, \pi]$, cf. \cite{Li2014,Nguyen2020}, as
\begin{align*}
	\bm{P}(\theta) = \begin{pmatrix} % bmatrix for square brackets
		\cos^2(\theta) 					& \sin^2(\theta) 					& - 2 \cos(\theta) \sin(\theta) 		\\
		\sin^2(\theta) 					 & \cos^2(\theta) 					& 2 \cos(\theta) \sin(\theta)			 \\
		 \cos(\theta) \sin(\theta) & -\cos(\theta) \sin(\theta)  & \cos^2(\theta) - \sin^2(\theta)
	\end{pmatrix}.
\end{align*}
The strong form of linear elasticity can be determined from $E_{bulk}$ in absence of volume loads and tractions and reads (cf., e.g., \cite{Wriggers2008})
%\begin{subequations}
\begin{align}
	\label{eqn:LinearElasticityStrongForm}
	\begin{aligned}
	\Div(\stress(\displacement)) &= \bm{0} \text{ in } \domain, \\
	\displacement &= \bar{\displacement} \text{ on } \partial \domain_D, \\
	\stress(\displacement) \bm{n} &= \bm{0} \text{ on } \partial \domain \setminus \partial \domain_D.
	\end{aligned}
\end{align}
%\end{subequations}
The part of the boundary $\partial \domain_D \subset \partial \domain$, $\partial \domain_D \neq \emptyset$ describes the part where \textsc{Dirichlet} boundary conditions are imposed\footnote{In the numerical experiments in \Cref{sec:NumericalExamples} \textsc{Dirichlet} boundary conditions are only imposed in $x$ or in $y$ direction to avoid the introduction of additional stresses. Only one point of the domain is constrained in both $x$ and $y$ direction. To avoid an overcomplicated and confusing notation we use the description in \eqref{eqn:LinearElasticityStrongForm} in any case.}.

In order to model the fracture behavior, \textsc{Griffith} has described a failure criterion based on the fracture toughness or \textsc{Griffith}'s criterion $G_c>0$ in \cite{Griffith1921}.
In the following, it is assumed that the crack surface $\crack$ is smooth enough
such that 
the fracture energy of $\crack$ can be described by
\begin{equation}
	\label{eqn:FractureEnergy}
E_{frac} = \int_\crack G_c \dd s
\end{equation}
as proposed in \cite{Francfort1998}.
 Therefore, the total energy of the system is given by
\begin{align*}
	E_{total} = E_{bulk} + E_{frac}.
\end{align*}
In view of the shape optimization approach, in some parts of the manuscript the total energy $E_{total}$ is also denoted by $J$.

%-----------------------------------------------------------------------
\subsection{Adaption of the fracture model to the phase-field setting}
\label{subsec:AdaptionFractureModelPhaseFieldSetting}
%----------------------------------------------------------------------- 
For the phase-field problem, a sharp-crack surface topology $\crack$ is regularized by the crack surface functional $\crack_{l_s}$ using scalar field variable $d \colon \mathcal{B} \to [0,1]$, $\bm{x} \mapsto d(\bm{x})$, which is introduced to indicate a crack if $d(\cdot) = 1$ and solid (intact) material if $d(\cdot) = 0$ along with the smooth transition among them, cf.~\cite{miehe+welschinger+hofacker10}.
The crack surface functional $\crack_{l_s}$ can be written as
\begin{align*}
\begin{split}
\crack_{l_s}(d) = \int_{\mathcal{B}} \gamma_l(d, \nabla d) \, \textrm{d}V %\\
\quad \mbox{with} \quad
\gamma_l(d, \nabla d) =  
\dfrac{1}{2l_s} d^2 + \dfrac{l_s}{2} \vert \vert\nabla d \vert\vert^2 \;,
\end{split}
%\label{s2-gamma_l}
\end{align*}
where $\gamma_l$ is the crack surface density function per unit volume of the solid and $l_s>0$ is the length scale parameter that governs the regularization.
The free energy required to generate the crack is taken from the classical definition using \textsc{Griffith}'s failure criteria as
\begin{align*}
	\begin{split}
		{\Psi}_{frac} = G_{c} \;\left( \dfrac{1}{2l_s} d^2 + \dfrac{l_s}{2} \vert \nabla d \vert^2 \right).
	\end{split}
\end{align*}
To describe a purely geometric approach to phase-field fracture, the regularized crack phase-field $d$ is obtained by a minimization principle of diffusive crack topology
\begin{align*}
d = \underset{d \in \mathcal{D}_\crack}{\arg\inf}\,
\crack_{l_s}(d)
%\; \mbox{with} \quad
%d = 1 \; \mbox{on} \; \crack\subset \mathcal{B}\TS{\text{still unclear}}
%\label{min-d-geo}
\end{align*}
with $\mathcal{D}_\crack \coloneqq \left\{ d \colon d(\crack) = 1 \right\}$. The \textsc{Euler} equations of the above variational principle are
\begin{subequations}
\begin{align*}
    d - l_s^2 \Delta d &= 0 \text{ in } \mathcal{B}, \\
    d &= 1 \text{ on } \crack, \\
    \nabla d \cdot \bm{n} &= 0 \text{ on } \partial \mathcal{B} \setminus \crack,
\end{align*}
\end{subequations}
where $\bm{n}$ is the outward normal on $\partial \mathcal{B}$, cf. \cite{Miehe2010+Hofacker+Welschinger}.
The local equation for the evolution of the crack phase-field in the domain $\mathcal{B}$ can be derived assuming global irreversibility condition of the crack evolution (cf. \cite{Miehe2010+Hofacker+Welschinger}) as
\begin{align*}
%\boxed{
	\frac{G_{c}}{l_s} \; ( d - l_s^2 \Delta d ) - 2 (1 - d) \Psi_0^+(\strain(\displacement)) = 0,
% }
%\TS{\text{still unclear};\text{why square brackets?}}
%\label{euler-eq-d}
\end{align*}
where the effective crack-driving force $ \Psi_0^+(\strain(\displacement))$ is based on strain splitting\footnote{For isotropic materials, the term $\Psi_0^+(\strain(\displacement))$ is defined as $\frac{K}{2} \left< \trace(\strain(\displacement)) \right>_+^2 + \mu \; (\strain^{dev}(\displacement) : \strain^{dev}(\displacement))$, where $K>0$ and $\mu>0$ denote bulk modulus and shear modulus, respectively, $\left< \cdot \right>_+$ is the \textsc{Macaulay} bracket, and $\strain^{dev}(\displacement)$ represents the deviatoric part of the strain tensor. For anisotropy, this is not usable---therefore, strain decomposition is not used.}.
%However, it has been shown in \cite{Lorenzis} that this strain decomposition can be omitted. Therefore, we choose
%\begin{align}
%	\mathcal{H} =  \max_{\displacement \in H^1(\domain, \R^2)} {\Psi_0}(\strain(\displacement)), %\; \ge 0
%	\label{driving-force}
%\end{align}
%which accounts for the irreversibility of the phase-field evolution by filtering out a maximum value of the crack driving state function $D$ in the full process history $s \in [0, t]$. \KW{We have already $s$ later for the piecwise smooth shapes.}
%The first term in \Cref{euler-eq-d} represents geometric resistance, the second term $ 2 (1 - d) \mathcal{H}$ is the driving force.
%Here, $\eta \ge 0$ is a parameter that characterizes the artificial/numerical viscosity of the crack propagation.
%is introduced that accounts for the irreversibility of the phase-field evolution by filtering out a maximum value of the crack driving state function $D$ in the full process history $s \in [0, t]$. 
%Herein, only the tensile part of the elastic energy $\psi_\textrm{el}^{+}$ is considered for the crack driving force, and it is given by 
%\begin{align}
%    \widetilde{\psi}_\textrm{el}^{+} = \frac{1}{2}\; K \left< \textrm{tr}(\boldsymbol{\epsilon})\right>_+^2 + \mu \; (\Tilde{\boldsymbol{\epsilon}} : \Tilde{\boldsymbol{\epsilon}})
%\end{align}
%The effective driving force $\mathcal{H}$ is updated at each \textsc{Gauss} point for every load step.
However, in line with the hybrid phase-field formulation from~\cite{ambati15} the total energy is degraded fully %., i.e., $\stress = g(d) \frac{\partial \Psi_0}{\partial \strain}$
in order to keep the computational cost comparable to the isotropic case. The degradation function is given as
\begin{align*}
	g(d) = (1 - d)^2 + k ,
\end{align*}
where $k>0$ is an algorithmic constant. Similar to \cite{KandWebHoech:2022} we choose $k = 10^{-7}$ for this paper. The parameter $k$ circumvents the full degradation of the free energy and leaves an artificial elastic rest energy density of $k \Psi_0(\strain(\displacement))$ when fully broken.
%-----------------------------------------------------------------------
\subsection{Adaption of the fracture model to the shape optimization setting}
\label{subsec:AdaptionFractureModelShapeOptimizationSetting}
%-----------------------------------------------------------------------
In contrast to the phase-field approach outlined in \Cref{subsec:AdaptionFractureModelPhaseFieldSetting}, an alternative approach to simulate fracture propagation has been proposed in \cite{Suchan2023} and is sketched in what follows. It is also based on the energy formulation but not using the phase-field $d$. 
Here, the problem is interpreted as a shape optimization problem with the goal of minimizing the energy~$E_{total}$ %in the system, 
constrained by a partial differential equation for the displacement $\displacement = (w_x, w_y)^\top \in H^1(\domain, \R^2)$. In our setting, the minimization of the energy is constrained by linear elasticity. The corresponding problem can be formulated as
\begin{subequations}
\begin{align}
	\label{eqn:ShapeOptimizationProblem1}
	\min_{\crack}& \int_\domain \frac{1}{2} \strain^V(\displacement) \cdot \bm{C} \cdot \strain^V(\displacement) \dd \bm{x} + \int_\crack G_c \dd s + \nu \int_\domain 1 \dd \bm{x} \\
	\label{eqn:ShapeOptimizationProblem2}
	&\text{s.t. } \int_\domain (\bm{C} \cdot \strain^V(\displacement)) \cdot \strain^V(\widetilde{\displacement}) \dd \bm{x} = 0 \quad \forall \widetilde{\displacement} \in H^1(\domain, \R^2) \\
	\label{eqn:ShapeOptimizationProblem3}
	&\text{ + appropriate boundary conditions.}
\end{align}
\end{subequations}
The last term in \eqref{eqn:ShapeOptimizationProblem1} describes a volume regularization term with $\nu > 0$ in order to incentivize that the fracture remains approximately one-dimensional.
We embed the optimization into a theoretical framework in order to benefit from optimization algorithms with established convergence properties. 

\begin{figure}[tbp]
	\centering
	\setlength\figureheight{5.5cm} 
	\setlength\figurewidth{5.5cm}
	\setlength\svgwidth{\figurewidth}
	%% Creator: Inkscape 1.0.2 (e86c8708, 2021-01-15), www.inkscape.org
%% Accompanies image file 'SENT_infinites_slit_svg-tex.pdf' (pdf, eps, ps)
\begingroup%
  \makeatletter%
  \providecommand\color[2][]{%
    \errmessage{(Inkscape) Color is used for the text in Inkscape, but the package 'color.sty' is not loaded}%
    \renewcommand\color[2][]{}%
  }%
  \providecommand\transparent[1]{%
    \errmessage{(Inkscape) Transparency is used (non-zero) for the text in Inkscape, but the package 'transparent.sty' is not loaded}%
    \renewcommand\transparent[1]{}%
  }%
  \providecommand\rotatebox[2]{#2}%
  \newcommand*\fsize{\dimexpr\f@size pt\relax}%
  \newcommand*\lineheight[1]{\fontsize{\fsize}{#1\fsize}\selectfont}%
  \ifx\svgwidth\undefined%
    \setlength{\unitlength}{319.26734819bp}%
    \ifx\svgscale\undefined%
      \relax%
    \else%
      \setlength{\unitlength}{\unitlength * \real{\svgscale}}%
    \fi%
  \else%
    \setlength{\unitlength}{\svgwidth}%
  \fi%
  \global\let\svgwidth\undefined%
  \global\let\svgscale\undefined%
  \makeatother%
  \begin{picture}(1,0.99791754)%
    \lineheight{1}%
    \setlength\tabcolsep{0pt}%
    \put(0.719732,0.81476681){\color[rgb]{0,0,0}\makebox(0,0)[t]{\lineheight{1.25}\smash{\begin{tabular}[t]{c}$\domain$\end{tabular}}}}%
    \put(0,0){\includegraphics[width=\unitlength,page=1]{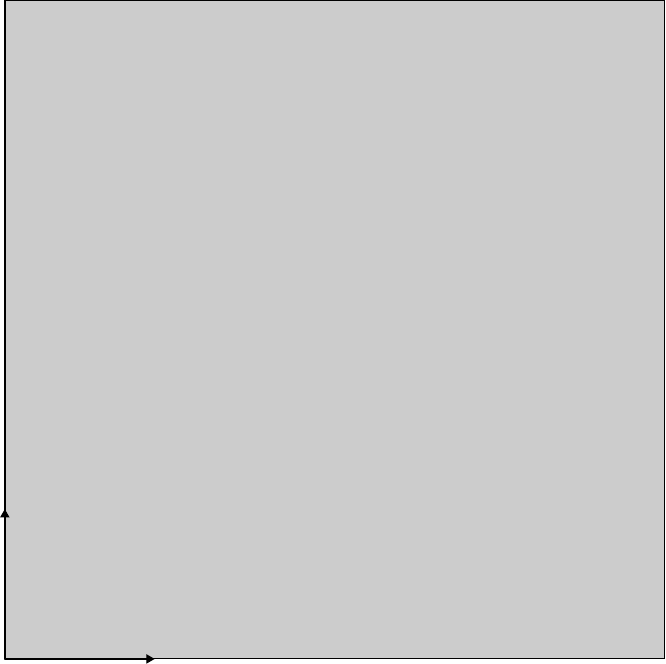}}%
    \put(0.49817388,0.54144999){\color[rgb]{0,0,0}\makebox(0,0)[t]{\lineheight{1.25}\smash{\begin{tabular}[t]{c}$\color{red}\crack$\end{tabular}}}}%
    \put(0,0){\includegraphics[width=\unitlength,page=2]{SENT_infinites_slit_svg-tex.pdf}}%
  \end{picture}%
\endgroup%
\quad%
	\newlength{\svgwidth}%
	\setlength\svgwidth{\figurewidth}%
	%% Creator: Inkscape 1.0.2 (e86c8708, 2021-01-15), www.inkscape.org
%% Accompanies image file 'SENT_triangle_slit_svg-tex.pdf' (pdf, eps, ps)
\begingroup%
  \makeatletter%
  \providecommand\color[2][]{%
    \errmessage{(Inkscape) Color is used for the text in Inkscape, but the package 'color.sty' is not loaded}%
    \renewcommand\color[2][]{}%
  }%
  \providecommand\transparent[1]{%
    \errmessage{(Inkscape) Transparency is used (non-zero) for the text in Inkscape, but the package 'transparent.sty' is not loaded}%
    \renewcommand\transparent[1]{}%
  }%
  \providecommand\rotatebox[2]{#2}%
  \newcommand*\fsize{\dimexpr\f@size pt\relax}%
  \newcommand*\lineheight[1]{\fontsize{\fsize}{#1\fsize}\selectfont}%
  \ifx\svgwidth\undefined%
    \setlength{\unitlength}{322.93002303bp}%
    \ifx\svgscale\undefined%
      \relax%
    \else%
      \setlength{\unitlength}{\unitlength * \real{\svgscale}}%
    \fi%
  \else%
    \setlength{\unitlength}{\svgwidth}%
  \fi%
  \global\let\svgwidth\undefined%
  \global\let\svgscale\undefined%
  \makeatother%
  \begin{picture}(1,0.98659915)%
    \lineheight{1}%
    \setlength\tabcolsep{0pt}%
    \put(0.72291081,0.80552572){\color[rgb]{0,0,0}\makebox(0,0)[t]{\lineheight{1.25}\smash{\begin{tabular}[t]{c}$\domain$\end{tabular}}}}%
    \put(0,0){\includegraphics[width=\unitlength,page=1]{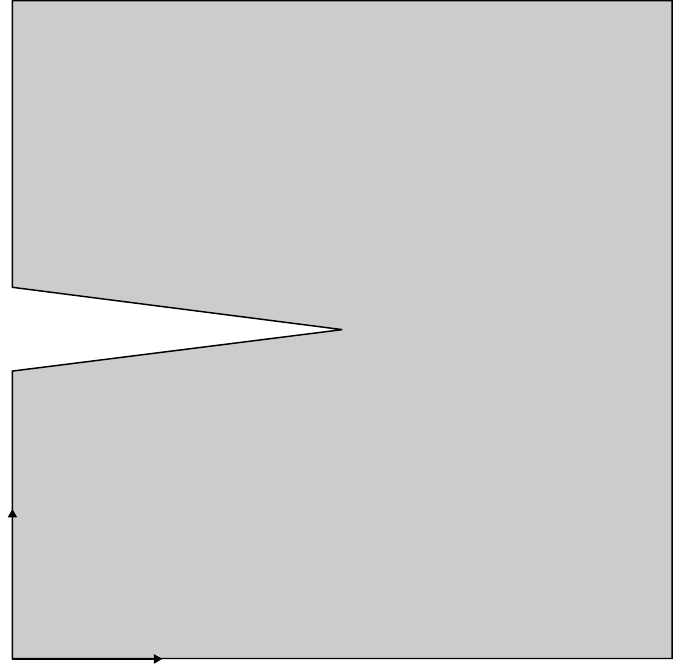}}%
    \put(0.0391694,0.58757494){\color[rgb]{0,0,0}\makebox(0,0)[lt]{\lineheight{1.25}\smash{\begin{tabular}[t]{l}$\color{blue}\bm{x}^0$\end{tabular}}}}%
    \put(0.5038656,0.53530886){\color[rgb]{0,0,0}\makebox(0,0)[t]{\lineheight{1.25}\smash{\begin{tabular}[t]{c}$\color{red}\crack$\end{tabular}}}}%
    \put(0,0){\includegraphics[width=\unitlength,page=2]{SENT_triangle_slit_svg-tex.pdf}}%
    \put(0.0391694,0.37079091){\color[rgb]{0,0,0}\makebox(0,0)[lt]{\lineheight{1.25}\smash{\begin{tabular}[t]{l}$\color{blue}\bm{x}^1$\end{tabular}}}}%
  \end{picture}%
\endgroup%
%
%	\includesvg[width=\figurewidth]{SENT_infinites_slit}\quad%
%%	\includegraphics{sent/mesh_iter0_compiled}
%	\includesvg[width=\figurewidth]{SENT_triangle_slit}
	\caption{Replacement of an infinitesimally thin fracture for the single-edge notch test (left) by an open curve as part of the boundary of $\domain$ (right).}
	\label{fig:ReplacementInfinitesimalFractureByTriangle}%
\end{figure}

In  \eqref{eqn:ShapeOptimizationProblem1}, the shape $\crack$ representing the fracture  is not specified so far.
Multiple options are available for the description of a shape---among these are landmark vectors \cite{Kendall1984}, plane curves \cite{Michor2006,Michor2007,Mio2006}, boundary contours of objects \cite{Ling2007,Rumpf2009}, multiphase objects \cite{Wirth2010}, characteristic functions of measurable sets \cite{Zolesio2007} or morphologies of images \cite{Droske2007}. We focus on the description of plane curves in this publication. For this, it has been motivated in \cite{Suchan2023} to replace the fracture by a curve towards the fracture tip and a curve back to the outer boundary in case the fracture has zero thickness. This setting is encountered, e.g., for the single-edge notch test, and is sketched in \Cref{fig:ReplacementInfinitesimalFractureByTriangle}. For an easier comparison we also use this mesh for the phase-field simulations.
Ideally, the space of all admissible shapes, the so-called shape space, would allow a vector space structure; this would make the formulation of algorithms relatively simple. However, many commonly used shape spaces do not provide this---among them is both, the shape space $B_e([0,1], \R^2)$ (cf. \cite{Michor2006}) containing all smooth embeddings of the unit interval into $\R^2$ excluding reparametrizations, as well as the shape space $M_s(B_e([0,1], \R^2)^N)$ from \cite{Pryymak2023}, where $s$ denotes the overall number of cracks. 
Both spaces admit an infinite dimensional \textsc{Riemann}ian manifold structure, which can then be exploited to formulate  gradient descent algorithms. An overview for shape optimization based on optimization on \textsc{Riemann}ian manifolds with a focus on engineering applications can be found in \cite{Radtke2023}. This concept has been extended in \cite{Geiersbach2023}, where the space $M_s^c$ has been introduced that can be used to describe multiple piecewise smooth closed curves. In the present publication, we require a shape space that allows multiple piecewise smooth open curves. Since the shape space $B_e([0,1], \R^2)$ does not allow the development of kinks in the shapes we adapt the definition of $M_s^c$ to obtain a space that allows $s \in \mathbb{N}$ piecewise smooth open curves\footnote{We use the superscript ``$o$'' for \emph{open} in contrast to ``$c$'' for \emph{closed} curves in this paper.} starting at point $\bm{x}_i^0$ and ending at point~$\bm{x}_i^1$, where $i$ indicates the respective crack:
\begin{align}
	\label{eqn:Definition_Mso}
\begin{aligned}
	M_s^o \coloneqq  
		& \left\{ \bm{\crack} = \left(\crack_1, \ldots, \crack_s \right) \in M_s(B_e([0,1], \R^2)^N)  \colon \right.\\
		& \left.\quad \crack_i\colon [0, 1) \to \R^2 \text{ injective with } \crack_{i}(0) = \widetilde{\crack}_{k_i}(0) = \bm{x}_i^0, \right. \\
		& \left.\quad \crack_{i}(1) = \widetilde{\crack}_{k_i+n_i-1}(1)  = \bm{x}_i^1    \text{ and } \widetilde{\crack}_{k_i+h}(1) = \widetilde{\crack}_{k_i+h+1}(0) \ \forall h = 0, \ldots, n_i-2 \right\},
		\end{aligned}
\end{align}
where
\begin{align*}
	\begin{aligned}
		%\begin{split}
		&M_s(B_e([0,1], \R^2)^N)\\
		&:= \left\{  \bm{\crack} = \left(\crack_1, \ldots, \crack_s \right) \colon \crack_i \in \prod_{l=k_i}^{k_i+n_i-1} B_e([0,1], \R^2) \ \forall i=1,\ldots,s, \right.\\
		& \hphantom{   := \,\,\,\, \,\,\, \bm{\crack} = \left(\crack_1, \ldots, \crack_s \right) \colon}\left.  \sum_{i=1}^s n_i = N, k_1=1, k_{i+1} = k_i + n_i \ \forall i=1, \ldots, s-1     \right\}
		%\end{split}
	\end{aligned}
\end{align*}
is the $s$-dimensional shape space on $B_e([0,1], \R^2)^N$ from \cite{Pryymak2023}. For all $i=1,\dots, s$ the piecewise smooth curve~$\crack_i$ in \eqref{eqn:Definition_Mso} consists of $n_i$ smooth curves~$\widetilde{\crack}_{k_i+h}$, $h=0,\ldots,n_i-2$. The start point of one smooth curve in \eqref{eqn:Definition_Mso} is connected to the end point of the previous smooth curve, i.e., the curve $\widetilde{\crack}_{k_i+h+1}$ starts at the point where the curve $\widetilde{\crack}_{k_i+h}$ ends. The start and end point of $\crack_i$ in \eqref{eqn:Definition_Mso} is defined as $\bm{x}_i^0$ and $\bm{x}_i^1$, respectively. 
% as
%\begin{align}
%	\label{eqn:Definition_Mso}
%	M_s^o \coloneqq \left\{ \begin{aligned}
%		& \bm{\crack} = \left(\crack_1, \ldots, \crack_s \right) : \crack_i \in \prod_{l=k_i}^{k_i+n_i-1} B_e([0,1], \R^2) \ \forall i=1,\ldots,s, \\
%		& \hphantom{\bm{\crack} = \left(\crack_1, \ldots, \crack_s \right) :} \sum_{i=1}^s n_i = N, k_1=1, k_{i+1} = k_i + n_i \ \forall i=1, \ldots, s-1, \\
%		& \crack_i: [0, 1) \to \R^2 \text{ injective with } \crack_{i}(0) = \widetilde{\crack}_{k_i}(0) = \bm{x}_i^0, \crack_{i}(1) = \widetilde{\crack}_{k_i+n_i-1}(1)  = \bm{x}_i^1 \text{ and} \\
%		& \quad \widetilde{\crack}_{k_i+h}(1) = \widetilde{\crack}_{k_i+h+1}(0) \ \forall h = 1, \ldots, n_i-2.\TS{\text{@KW: should this be $h=0,\ldots,n_i-2$?}}
%	\end{aligned} \right\}
%\end{align}
A kink in the piecewise smooth curve $\crack_i$ can then develop at these connecting points.
In some applications, the fracture is defined by only one open curve. In this case, the definition \eqref{eqn:Definition_Mso} simplifies to
\begin{align*}
%	\label{eqn:Definition_Mso_1curve}
	\begin{aligned}
	M_1^o \coloneqq 
		&\left\{ \crack \in M_1(B_e([0,1], \R^2)^N) \colon \right.  \\
		& \left. \quad \crack\colon [0, 1) \to \R^2 \text{ injective with } \crack(0) = \widetilde{\crack}_{1}(0) = \bm{x}^0, \crack(1) = \widetilde{\crack}_{N}(1)  = \bm{x}^1 \text{ and} \right. \\
		& \left. \quad \widetilde{\crack}_{h}(1) = \widetilde{\crack}_{h+1}(0) \ \forall h = 1, \ldots, N-1 \right\}.
	\end{aligned}
\end{align*}
In order to establish a measure of distance between two different elements of $M_s^o$ and to determine how one element is changed to another, a \textsc{Riemann}ian metric on the shape space is required. The \textsc{Steklov}-\textsc{Poincaré} metric (cf. \cite{Schulz2016,Siebenborn2017}) is used in this publication because of its advantageous properties regarding the quality of the computational domain discretization. Due to the usage of the \textsc{Steklov}-\textsc{Poincaré} metric, we obtain a vector field $\bm{V} \colon \domain \to \R^2$ in an appropriate \textsc{Hilbert} space $\mathcal{V}$ by solving the so-called deformation equation
\begin{align}
	\label{eqn:DeformationEquationGeneral}
	a(\bm{V}, \bm{W}) = \mathrm{d} J(\crack)[\bm{W}] \quad \forall \bm{W} \in \mathcal{V}
\end{align}
with a symmetric and coercive bilinear form $a\colon \mathcal{V} \times \mathcal{V} \to \R$ and the shape derivative $\mathrm{d} J$ \cite[Eq.~(4.4)]{Schulz2016}.
Please note that $J$ equals the total energy $E_{total}$ according to~\Cref{subsec:FractureModel}.
The resulting vector field $\bm{V}$ can then be used to update the domain as
\begin{align*}
	\domain_\tau = \{ \bm{x} \in \domain \colon \bm{x} + \tau \bm{V}(\bm{x}) \}
\end{align*}
with $\tau \geq 0$. The update of the curve $\crack \subset \domain$ is, accordingly,
\begin{align*}
	\crack_\tau = \{ \bm{x} \in \crack \colon \bm{x} + \tau \!\left.\bm{V}(\bm{x})\right|_{\crack} \}.
\end{align*}
The shape derivative of $J$ at $\crack$ in direction $\bm{W}$ denoted by $\mathrm{d} J(\crack)[\bm{W}]$ describes the change in energy due to the change of $\crack$ (and $\domain$, respectively) and is defined as
\begin{align*}
	\mathrm{d} J(\crack)[\bm{W}] = \lim_{\tau \to 0^+} \frac{J(\crack_\tau) - J(\crack)}{\tau}
\end{align*}
if it exists and is linear and continuous in $\bm{W}$. As described in \cite{Suchan2023}, the shape derivative for our problem can be calculated with a material derivative approach \cite{Berggren2009} as
\begin{align}
%	\begin{aligned}
	\nonumber
	\mathrm{d} J(\crack)[\bm{W}] = &\int_\domain -\frac{1}{2} \left( \nabla \displacement \nabla \bm{W} + \left( \nabla \displacement \nabla \bm{W} \right)^\top \right) : \stress(\displacement) + \operatorname{div}(\bm{W}) \cdot \left( \frac{1}{2} \stress(\displacement) : \strain(\displacement) \right) \dd \bm{x} \\
	&+ \frac{1}{2} \int_\crack G_c \, \kappa \bm{W}^\top \bm{n} \dd s - \int_\crack \nu \bm{W}^\top \bm{n} \dd s
%	\end{aligned}
\end{align}
with $\kappa$ as the curvature of $\crack$.
The bilinear form $a(\bm{V}, \bm{W})$ in \eqref{eqn:DeformationEquationGeneral} can be chosen arbitrarily and \eqref{eqn:DeformationEquationGeneral} is solely used to obtain the vector field $\bm{V}$. Common choices include the \textsc{Poisson} equation or the bilinear form of linear elasticity. Similar to \cite{Suchan2023} we use the bilinear form of linear elasticity. Solving for $\bm{V}$ does not replace the linear elasticity equation to obtain $\displacement$ but instead is an additional step to obtain the gradient. The problem to obtain for $\bm{V} \in \mathcal{V}$ therefore reads
\begin{align}
	\label{eqn:DeformationEquationLinElas}
	\begin{aligned}
	\tunderbrace{\int_{\domain} \left(\widetilde{\mathbb{C}} : \strain(\bm{V}) \right) : \strain(\bm{W}) \dd \bm{x}}{a(\bm{V}, \bm{W})} &= \mathrm{d} J(\crack)[\bm{W}] \quad \forall \bm{W}\in \mathcal{V}.
	\end{aligned}
\end{align}
Following~\cite{Schulz2016,Siebenborn2017} we choose the second \textsc{Lamé} parameter in $\widetilde{\mathbb{C}}$ to obtain $\bm{V}$ as
%
%In \cite{Schulz2016,Siebenborn2017} an isotropic material stiffness tensor $\mathbb{C}$ is chosen in \eqref{eqn:DeformationEquationLinElas} and the \textsc{Lamé} parameters are determined by prescribing values on $\crack$ and on $\partial \domain \setminus \crack$ and then solving a Poisson problem to extend the parameters into the domain, which is also followed in this manuscript. The material parameters $\tilde{\lambda}$ and $\tilde{\mu}$ that describe $\mathbb{C}$ in \eqref{eqn:DeformationEquationLinElas} do not have to correspond to the material parameters $\lambda$ and $\mu$ of the physical system, and instead can be choosen freely. In fact, the material stiffness tensor of the physical system does not even have to model isotropic behavior. Further, it has been suggested in~\cite{Siebenborn2017} to choose a larger $\tilde{\mu}$ on $\crack$ which helps with the quality of the discretization of the domain during the optimization. Therefore, we choose 
$\tilde{\mu}=5$ on $\crack$ and $\tilde{\mu}=1$ on $\partial \domain \setminus \crack$ and perform an extension to obtain intermediate values inside of $\domain$. The first \textsc{Lamé} parameter is chosen as $\tilde{\lambda}=10$ everywhere. These values are unrelated to the material parameters in the balance of linear momentum to obtain $\displacement$.
Solving \eqref{eqn:DeformationEquationLinElas} further requires appropriate boundary conditions. Since only deformations of $\crack$ should affect the objective functional that is to be minimized, an appropriate choice is to not allow deformations of $\partial \domain \setminus \crack$, which yields a homogenous \textsc{Dirichlet} boundary condition on $\partial \domain \setminus \crack$. We further impose a homogenous \textsc{Dirichlet} boundary condition on the part of $\crack$ that fulfills $0<x<0.48$ since no fracture growth can be expected in this area, and a homogenous \textsc{Neumann} boundary condition on the rest of $\crack$, which we denote as $\hat{\crack}$. The previously mentioned \textsc{Hilbert} space $\mathcal{V}$ therefore can be chosen as
\begin{align*}
	\mathcal{V} \coloneqq \left\{ \bm{V} \in H^1(\domain, \R^2) \colon \bm{V}=\bm{0} \text{ on } \partial \domain \setminus \hat{\crack} \right\}.
\end{align*}

Using the approach described above the fracture irreversibility condition is not ensured. Thus, for low mechanical loads on the physical system the fracture would actually shrink, which is an unphysical behavior, since the fracture cannot heal itself after developing. Similar to \cite{Suchan2023}, we project the resulting vector field $\bm{V}$ from \eqref{eqn:DeformationEquationLinElas} onto the set
\begin{align}
	\label{eqn:IrreversibilityProjectionSet}
	\left\{ \bm{V} \in H^1(\domain, \R^2) \colon \bm{V}^\top \bm{N} \leq 0, \bm{V}=\bm{0} \text{ on } \partial \domain \setminus \hat{\crack} \right\}
\end{align}
to ensure irreversibility, where $\bm{N}$ describes the unit outward normal vector on $\crack$ extended into the domain $\domain$. This is performed by solving the Eikonal equation
\begin{align}
	\label{eqn:EikonalEquation}
	\varepsilon \Delta \Phi(\bm{x}) + \left| \nabla \Phi(\bm{x}) \right| = 1 \quad \text{in } \domain, \quad \Phi(\bm{x}) = 0 \quad \text{on } \crack, \quad \frac{\partial \Phi(\bm{x})}{\partial \bm{n}} = 0 \quad \text{on } \partial \domain \setminus \crack
\end{align}
stabilized by an $\varepsilon=2 \cdot 10^{-3}$. The normal vector extended into the domain $\domain$ can then be determined as $\bm{N} = \nabla \Phi(\bm{x})$. The irreversibility is enforced by replacing any $\bm{V} \in \mathcal{V}$ that violates $\bm{V}^\top \bm{N} \leq 0$ with $\bm{V}=\bm{0}$. This equates to a projected gradient descent algorithm.

%=======================================================================
\section{Numerical implementation}
\label{sec:Implementation}
%=======================================================================
In this section, numerical implementations are compared for both methods in order to get detailed insight about different approaches.
%-----------------------------------------------------------------------
\subsection{Implementation of the phase-field approach}
\label{subsec:ImplementationPhaseField}
%-----------------------------------------------------------------------
%\TS{todo: strong form of linear elasticity with phase-field approach}
As described in~\cite{Hirshikesh2019} the weak form of the bulk response in the \textsc{Euler}ian setting is given as
\begin{align}
	\label{eq:Weak_form_deformation}
	\int_{\domain} g(d) \, {\stress} : \nabla \delta \displacement 
	%- \bm{\gamma} \cdot \delta \displacement
	 \dd \bm{x} &= 0 \quad \forall \delta \displacement \in H^1(\domain, \R^2), \\
	%-\int_{\partial \domain } \bar{\bm{t}} \cdot \delta \displacement \: da = 0,
%\end{align}
%\begin{align}
	\label{eq:Weak_form_phase-field}
	\int_{\domain} \delta d \left( \dfrac{G_c}{l_s} + 2 \mathcal{H} \right) d - 2 \mathcal{H} \delta d + \nabla \delta d \cdot G_c l_s \nabla  d \dd \bm{x} &=  0 \quad \forall \delta d \in H^1(\domain, \R).
\end{align}
The finite element method is used to discretize the domain $\domain = \mathop{\bigcup}_{e=1}^{N_e} \; \domain^e$ with ${N}_e$ as the number of elements.
The residuum in the domain $\domain$ with the finite element discretization of the elements using linear shape functions with $\bm{N}_{(\cdot)}(\bm{x})$ as the shape functions and $\bm{B}_{(\cdot)}(\bm{x})$ as their respective derivatives, are given as 
\begin{align}
	\label{eq:Residual_deformation}
	\bm{r}_{\displacement}
	%(\bm{d}_{\displacement}, \bm{d}_{d}) 
	&=  \mathop{\bigcup}_{e=1}^{N_e} \int_{\domain^e} g(d) \, { \bm{B}_{\displacement}^\top \: \stress^V}
	%- \bm{N}_{\displacement}^T\: \bm{\gamma} 
	\dd \bm{x} = \bm{0}, \\
	%- \mathop{\bigcup}_{s=1}^{N_e} \int_{\partial \domain^e} \bm{N}_{\displacement}^T \: \bar{\bm{t}} \: da = 0,
%\end{align}
%\begin{align}
	\label{eq:Residual_phase-field}
	r_d 
	%(\bm{d}_{\displacement}, \bm{d}_{d}) 
	&= \mathop{\bigcup}_{e=1}^{N_e} \int_{\domain^e} \bm{N}_{d}^\top \left[\left( \dfrac{G_c}{l_s} + 2 \mathcal{H} \right) d - 2 \mathcal{H} \right] + \bm{B}_{d}^\top l_s^2\nabla d \dd \bm{x} = 0.
\end{align}
This multi-field problem is solved using \textsc{Newton-Raphson} iterations in a staggered manner as depicted in \Cref{alg:algorithmPhaseField}.
Therein, a sequence of two linear subproblems for the successive update of the phase-field fracture and the displacement field is shown. 
\begin{algorithm}
	\caption{Numerical implementation for phase-field approach}
	\label{alg:algorithmPhaseField}
	\begin{algorithmic}[1]
		\State \textbf{Input:} displacement increment $\triangle \bar{\displacement} > \bm{0}$, max. number of staggered iterations $k_{\max}>0$, staggered tolerance $tol_{stag} = 10^{-3}$
		%\State solution  $(\displacement_{n-1},d_{n-1})$ form step n-1.
		\State \textbf{Initialization:} $\displacement = \bm{0}$, $d=0$, $\bar{\displacement}=\bm{0}$
		\While{$\bar{\displacement} < \bar{\displacement}_{max}$}
		\State \textbf{Initialization:} $k=0$, $(\displacement^0, d^0)= (\displacement,d)$
%		\State Staggered iteration: for $k \geq 1$
		\While{$k < k_{max}$} 		
		\State Solve \eqref{eq:Residual_deformation} to obtain $\displacement^k \in H^1(\domain, \R^2)$ with $d^{k-1}$ and $\bar{\displacement}$
		\State Solve \eqref{eq:Residual_phase-field} to obtain $d^k \in H^1(\domain, \R)$ with $\displacement^k$
		\State Compute errornorm for obtained pair $(\displacement^{k}, d^k)$
		\State \algorithmicif\ errornorm $\leq$ $tol_{stag}$ \algorithmicthen\ $(\displacement, d) \leftarrow (\displacement^k,d^k$) and \algorithmicbreak
		\State $k \leftarrow k + 1$
		\EndWhile
		\State $\bar{\displacement} \leftarrow \bar{\displacement} + \triangle \bar{\displacement}$
		\EndWhile
		\State \textbf{Output:} solution $(\displacement,d)$.
	\end{algorithmic}
\end{algorithm}

\subsection{Implementation of the shape optimization approach}

Similar to the phase-field approach, however without the phase-field variable $d$, the weak form for the displacement is given by linear elasticity (cf., e.g., \cite{Bathe2007}) as
\begin{align*}
%	\label{eq:Weak_form_deformation_ShapeOpt}
	\int_{\domain} {\stress(\displacement)} : \nabla \delta \displacement  \dd \bm{x} = 0 \quad \forall \delta \displacement \in H^1(\domain, \R^2).
\end{align*}
The finite element method with the same initial discretization and shape functions as in \Cref{subsec:AdaptionFractureModelPhaseFieldSetting} is used to solve the PDE. Therefore, the residuum reads
\begin{align}
	\label{eq:Residual_deformation_ShapeOpt}
	\bm{r}_{\displacement}
	%(\boldsymbol{d}_{\displacement}, \boldsymbol{d}_{d}) 
	&=  \mathop{\bigcup}_{e=1}^{N_e} \int_{\domain^e} { \bm{B}_{\displacement}^\top \: \stress^V(\displacement)}
	%- \boldsymbol{N}_{\displacement}^T\: \boldsymbol{\gamma} 
	\dd \bm{x} = \bm{0}
\end{align}
and is solved by classical LU\footnote{`LU' is the abbreviation of `lower-upper.'} factorization. Then, an additional linear PDE is solved to obtain the vector field from the shape derivative to update the computational mesh as described in \eqref{eqn:DeformationEquationLinElas}. This yields the same weak form and discretization as for the displacement $\displacement$ but with different material parameters. The specific choice of material parameters is described in \Cref{sec:NumericalExamples}. The Eikonal equation \eqref{eqn:EikonalEquation} is also solved using the same discretization with the standard \textsc{Newton} method, from which the extension of the normal vector into the domain is obtained by projecting $\nabla \Phi$ to linear shape functions again. The projection onto the set~\eqref{eqn:IrreversibilityProjectionSet} is performed by computing the scalar product of $\bm{V}$ and $\bm{N}$ at each degree of freedom (DOF) and if the scalar product is larger than zero, then $\bm{V}=\bm{0}$ at this DOF. The projected vector field $\bm{V}$ is then applied to $\domain$ using a constant step size of $\tau=10^{-2}$. For our simulations, we choose a maximum number of iterations of $k_{\max}=5000$, however this limit was not reached in any simulation shown in \Cref{sec:NumericalExamples}. The whole process is sketched in \Cref{alg:algorithmShapeOpt}.

\begin{algorithm}
	\caption{Numerical implementation for shape optimization}
	\label{alg:algorithmShapeOpt}
	\begin{algorithmic}[1]
		\State \textbf{Input:} Initial domain $\domain^1 \subset \R^2$, step size $\tau>0$, max. number of iterations $k_{\max}>0$, displacement increment $\triangle \bar{\displacement} > \bm{0}$
		\State \textbf{Initialization:} $\bar{\displacement}=\bm{0}$
		\While{$\bar{\displacement} < \bar{\displacement}_{\max}$}
		\State \textbf{Initialization:} $k\coloneqq1$
		\While{$k < k_{\max}$}
		\State Solve \eqref{eq:Residual_deformation_ShapeOpt} to obtain $\displacement \in H^1(\domain, \R^2)$ with $\bar{\displacement}$
		\State Compute $J(\domain^k)$
		\State Generate vector field $\bm{V}$ by solving \eqref{eqn:DeformationEquationLinElas}
		\State Project $\bm{V}$ onto \eqref{eqn:IrreversibilityProjectionSet}
		\State $\domain^{k+1} \leftarrow \{ \bm{x} \in \domain^k \colon \bm{x} + \tau \bm{V}(\bm{x}) \}$
		\State \algorithmicif\ {$J(\domain^{k+1}) \geq J(\domain^k)$}\ \algorithmicthen\ \algorithmicbreak
		\State $k\leftarrow k+1$
		\EndWhile
		\State $\bar{\displacement} \leftarrow \bar{\displacement} + \triangle \bar{\displacement}$
		\EndWhile
	\end{algorithmic}
\end{algorithm}

%=======================================================================
\section{Numerical examples}
\label{sec:NumericalExamples}
%=======================================================================
%\KW{What do you mean by 'representative'? Representative for what? Maybe only 'Numerical examples'?}

In this section, we show the numerical results from the Algorithms \ref{alg:algorithmPhaseField} and \ref{alg:algorithmShapeOpt} for two different boundary value problems. In \Cref{subsec:BVP1}, we investigate the V-notch model as described in~\cite{Boyce2016} under tension and compare boundary forces and fracture path between the two different approaches. Then, \Cref{subsec:BVP2} demonstrates the numerical results for fracture propagation of the single-edge notch tension test with anisotropic material behavior.

%-----------------------------------------------------------------------
%\subsection{BVP1: I-beam - Isotropy}
%\label{subsec:BVP1}
%%-----------------------------------------------------------------------
%
%\begin{figure}[tbp]
%    \centering
%    \includegraphics[height=0.5\linewidth]{Figs/Ibeam_meshR2.png}
%    \caption{Mesh for I beam: Nodes - 2086, Elements - 4188}
%    \label{fig:BVP1_Ibeam_Mesh_R1}
%\end{figure}

%-----------------------------------------------------------------------

\subsection{Boundary Value Problem 1: V-notch -- Isotropy}

\label{subsec:BVP1}
%-----------------------------------------------------------------------
\begin{figure}[h]
	\centering
	\setlength\figureheight{7cm} 
	\setlength\figurewidth{.95\textwidth}
	\includegraphics{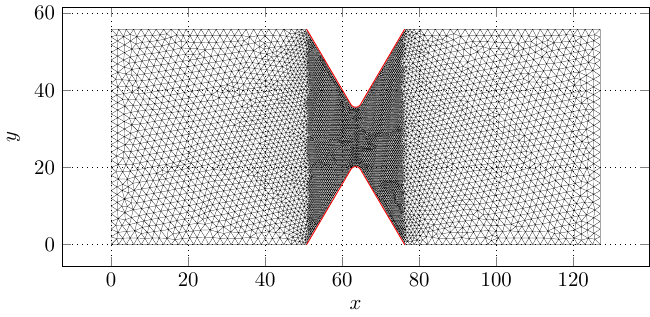}
%		\includetikz{vnotch/mesh_iter0}
	\caption{Mesh for V-notch: Nodes - 6900, Elements - 13434}
	\label{fig:BVP1_Vnotch_Mesh}%
\end{figure}
For this numerical experiment, we use the double V-notch plate with the geometric properties and computational mesh shown in~\Cref{fig:BVP1_Vnotch_Mesh}. The discretization of the model contains $13\,434$~triangular elements. The mesh is refined in the central region where the crack is expected to develop. The red lines denote the two curves $\crack_1$ and $\crack_2$, cf.~\eqref{eqn:Definition_Mso}. The notch is at an angle of~$60^\circ$ and has a rounded tip with a radius of $1.3\,\text{mm}$.
The left edge is fixed in $y$-direction and the top left corner node is fixed in both $x$ and $y$ direction. Displacement control on the right edge is used for this simulation with a displacement increment of $\triangle \bar{\displacement} = 10^{-4}\,\text{mm}$ up to failure for phase-field simulations and $\triangle \bar{\displacement} = 10^{-2}\,\text{mm}$ for shape optimization simulations.
The material parameters are chosen as $\lambda = 121.15\,\text{GPa}$, $\mu = 80.77\,\text{GPa}$, $G_c = 2.7\,\text{N/mm}$. For the phase-field method, the length scale parameter was chosen as $l_s = 0.01\,\text{mm}$.

\begin{figure}[tbp]
	\centering
	\setlength\figurewidth{.45\textwidth}%
	\includegraphics[width=\figurewidth]{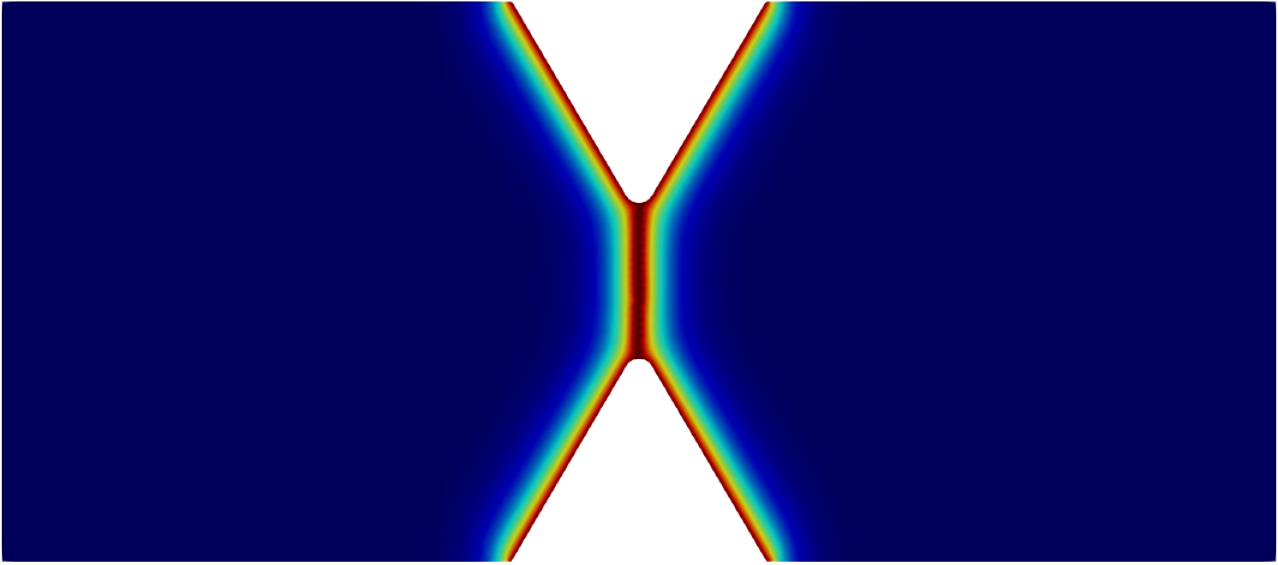}\qquad%
	\includegraphics{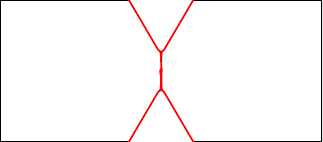}%
%	\includetikz{vnotch/crackpaths/bounds_iter8}%
	\caption{Fracture evolution results: phase-field (left) and shape optimization (right).}
	\label{fig:BVP1_Vnotch_Results}%
\end{figure}

\Cref{fig:BVP1_Vnotch_Results} shows the comparison of the fracture paths between the phase-field and the shape optimization approach.
On the left-hand side of the figure, blue represents the intact material and red represents the fractured regions. On the right-hand side, the curves representing the fracture are shown in red.
From the results it can be seen that the crack paths from both approaches are in good agreement.

%\TS{TODO: dont forget to describe start and end point of curve}

%-----------------------------------------------------------------------
\subsection{Boundary Value Problem 2: Single-edge notch test -- Anisotropy}
\label{subsec:BVP2}

%\cite[Sec. 3.5]{Boresi1993}

\begin{figure}[h]
	\centering
	\setlength\figureheight{10cm} 
	\setlength\figurewidth{10cm}
	\includegraphics{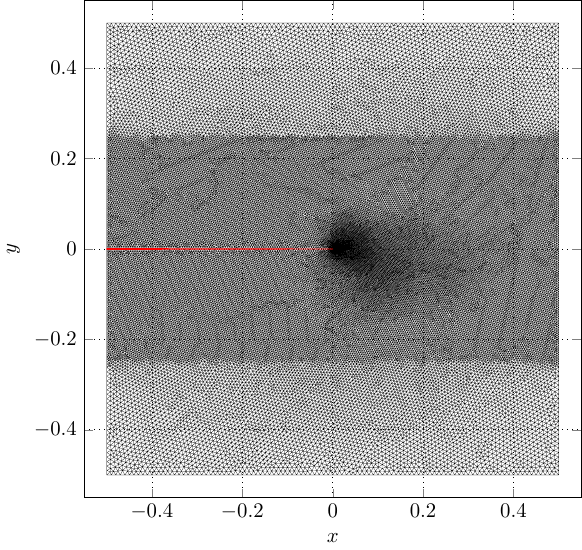}
	%	\includetikz{sent/mesh_iter0}
	\caption{Mesh for SENT: Nodes - 34863, Elements - 69014}
	\label{fig:BVP2_SENT_Mesh}%
\end{figure}

We now model the single-edge notch test (SENT) for testing and validating the extension of the shape opimization approach towards anisotropic elastic media, e.g., monocrystals with a predefined orientation angle. 
For this numerical experiment we used a square plate with a width and height of $1\,\text{mm}$ and a pre-existing crack from the left edge until the middle of the plate. The computational mesh consists of $69\,014$ elements and it is refined in the middle section where the crack is expected to develop, as shown in \Cref{fig:BVP2_SENT_Mesh}. Furthermore, the curve representing the crack is again shown in red and has a width of $2 \cdot 10^{-3}\,\text{mm}$ on the left edge of the domain (cf. \Cref{fig:ReplacementInfinitesimalFractureByTriangle}).

The elastic material stiffness matrix $\bm{C}^{ref}$ (in \textsc{Voigt} notation) for the reference coordinate system is defined as
\begin{align*}
	 \bm{C}^{ref} = 
	\begin{pmatrix}
		65 & 20 & 0 \\
		20 & 260 & 0 \\
		0 & 0 & 30 
	\end{pmatrix} \: \text{GPa}.
\end{align*}
Following~\cite{Nguyen2020}, the critical energy is considered independent of the orientation angle and is taken as $G_c = 1\,\text{N/mm}$. 
For the phase-field simulations, the legth scale parameter is chosen as $l_s = 0.005\,\text{mm}$. 
Simulations are performed using displacement control with an increment of $\triangle \bar{\displacement} = 2 \cdot 10^{-5}\,\text{mm}$ on the top edge when using the phase-field and $\triangle \bar{\displacement} = 10^{-4}\,\text{mm}$ when using the shape optimization approach. Further, we impose a homogenous \textsc{Dirichlet} boundary condition in $x$ direction on all nodes on the bottom edge and fix the bottom left corner of the domain in both $x$ and $y$ direction.
Numerical computations are performed for brittle fracture with different anisotropy angles $0^\circ$, $30^\circ$, $60^\circ$ and $90^\circ$ measured from the horizontal axis.

\begin{figure}[htb]
	\centering
	\setlength\figureheight{.3\textwidth}  % könnten auch global gesetzt werden.
	\setlength\figurewidth{.3\textwidth}  % könnten auch global gesetzt werden.
	\begin{subfigure}[T]{\textwidth}%
		\centering%
		\includegraphics[width=.87\figurewidth]{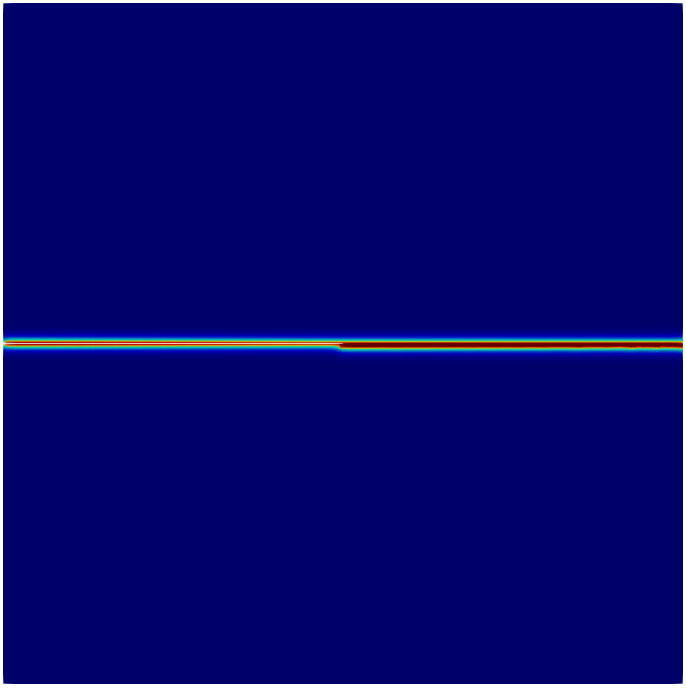}\qquad%
		\includegraphics{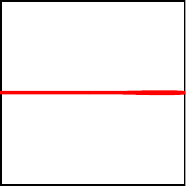}%
%		\includetikz{sent/fractures_shapeOpt/0degrees}%
		\vspace*{2pt}
		\caption{$\theta=0^\circ$}%
		\label{fig:FractureShapeOpt_0degrees}%
	\end{subfigure}\\[10pt]%
	\begin{subfigure}[T]{\textwidth}%
		\centering%
		\includegraphics[width=.87\figurewidth]{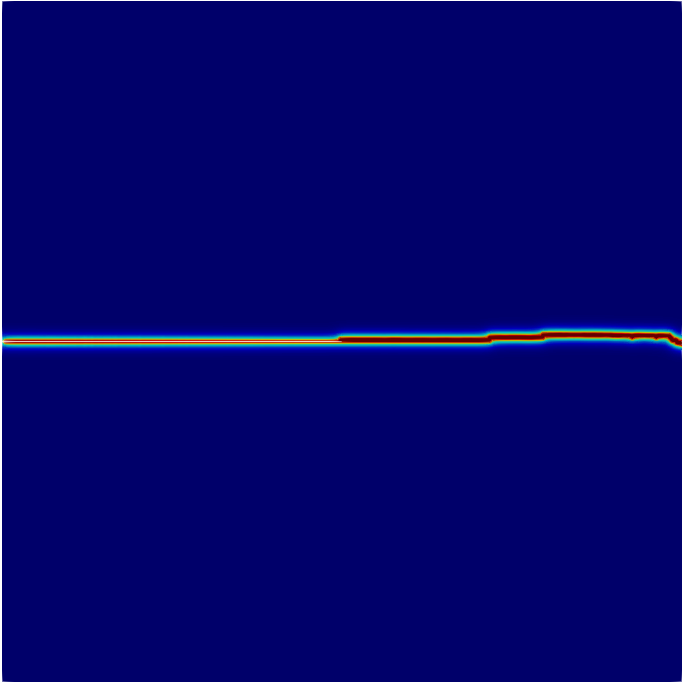}\qquad%
		\includegraphics{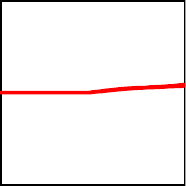}%
%		\includetikz{sent/fractures_shapeOpt/30degrees}%
		\vspace*{2pt}
		\caption{$\theta=30^\circ$}%
		\label{fig:FractureShapeOpt_30degrees}%
	\end{subfigure}\\[10pt]%
	\begin{subfigure}[T]{\textwidth}%
		\centering%
		\includegraphics[width=.87\figurewidth]{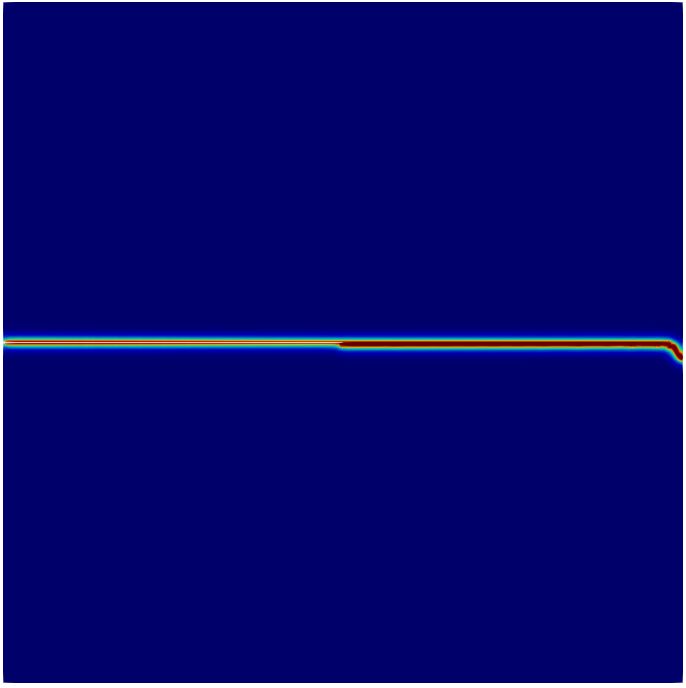}\qquad%
		\includegraphics{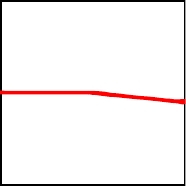}%
%		\includetikz{sent/fractures_shapeOpt/60degrees}%
		\vspace*{2pt}
		\caption{$\theta=60^\circ$}%
		\label{fig:FractureShapeOpt_60degrees}%
	\end{subfigure}\\[10pt]%
	\begin{subfigure}[T]{\textwidth}%
		\centering%
		\includegraphics[width=.87\figurewidth]{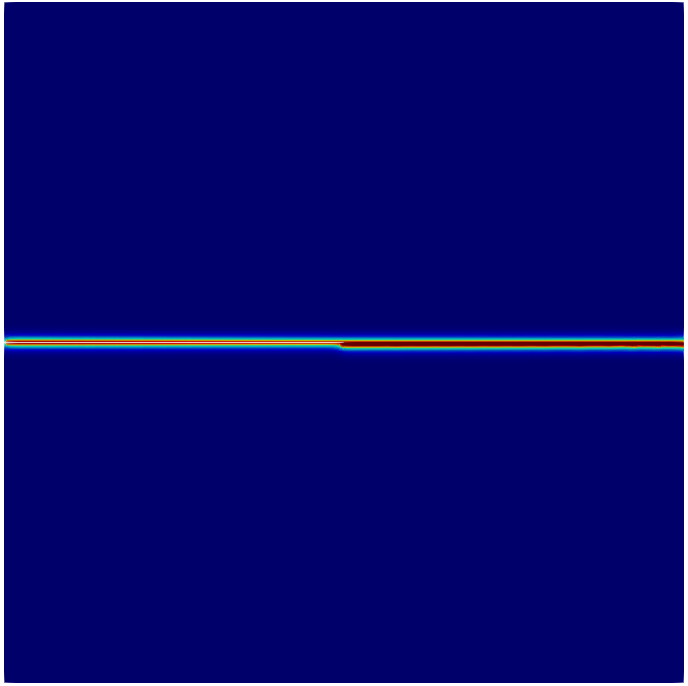}\qquad%
		\includegraphics{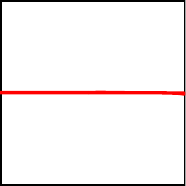}%
%		\includetikz{sent/fractures_shapeOpt/90degrees}%
		\vspace*{2pt}
		\caption{$\theta=90^\circ$}%
		\label{fig:FractureShapeOpt_90degrees}%
	\end{subfigure}
	\caption{Fracture paths for different angles $\theta$ of the anisotropic material.}
	\label{fig:FractureShapeOpt}
\end{figure}

\begin{figure}[tbp]
	\centering
	\setlength\figureheight{8cm} 
	\setlength\figurewidth{.95\textwidth}
	\includegraphics{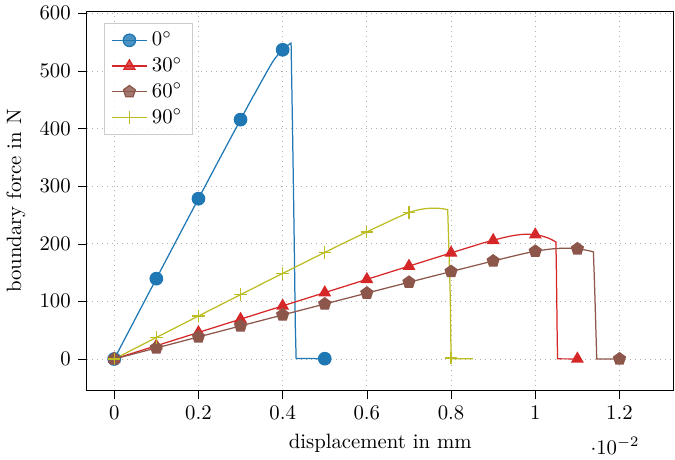}
%	\includetikz{sent/force_vs_displacement_all}
	\caption{Boundary force for SENT anisotropy with the phase-field approach.}
	\label{fig:BVP2_phasefield_boundaryForce}%
\end{figure}

\begin{figure}[tbp]
	\centering
	\setlength\figureheight{8cm} 
	\setlength\figurewidth{.95\textwidth}
	\includegraphics{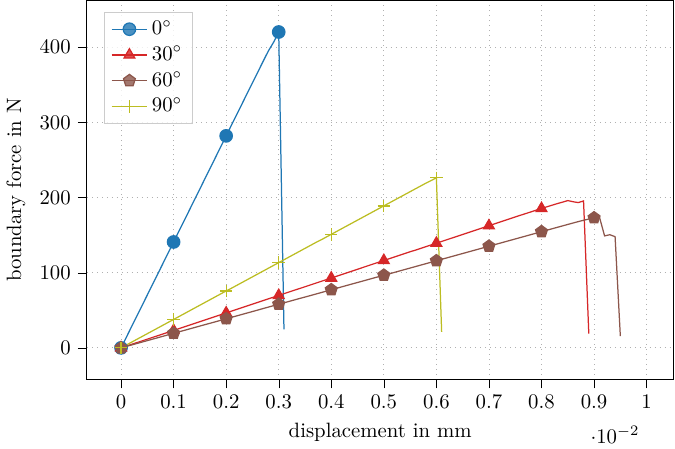}
%	\includetikz{sent/boundaryForceY_Increment1e-4}
	\caption{Boundary force for SENT anisotropy with the shape optimization approach.}
	\label{fig:BVP2_SENT_BoundaryForce_Increment1e-4}%
\end{figure}

\begin{figure}[tbp]
	\centering
	\setlength\figureheight{8cm} 
	\setlength\figurewidth{.95\textwidth}
		\includegraphics{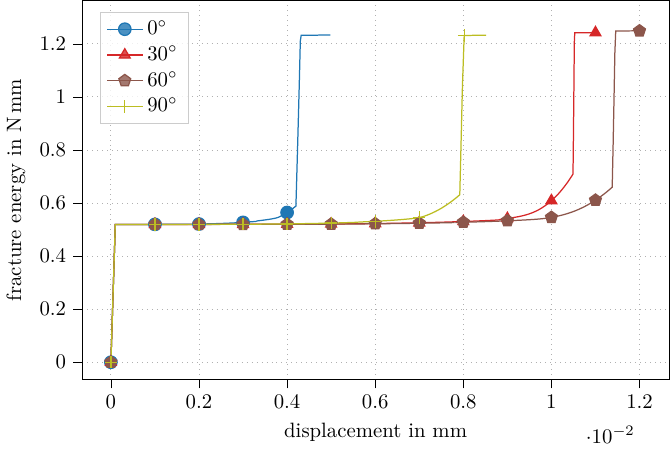}
%	\includetikz{sent/Fracture_Integral_all}
	\caption{Fracture energy for SENT anisotropy with the phase-field approach.}
	\label{fig:BVP2_phasefield_FractureEnergy}%
\end{figure}

\begin{figure}[tbp]
	\centering
	\setlength\figureheight{8cm} 
	\setlength\figurewidth{.95\textwidth}
		\includegraphics{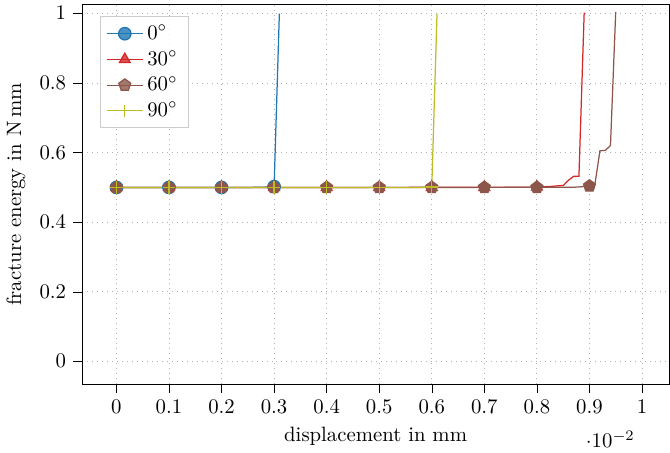}
%	\includetikz{sent/objective_functional_crack}
	\caption{Fracture energy for SENT anisotropy with the shape optimization approach.}
	\label{fig:BVP2_ShapeOpt_FractureEnergy}%
\end{figure}

The fracture paths for all anisotropy angles are depicted in~\Cref{fig:FractureShapeOpt}. 
We observe good qualitative agreement of the phase-field and the shape optimization approach, as the direction of fracture propagation between the two approaches is identical. The boundary forces for the phase-field simulations are shown in \Cref{fig:BVP2_phasefield_boundaryForce} and for the shape optimization simulations in \Cref{fig:BVP2_SENT_BoundaryForce_Increment1e-4}. Also here, we observe good qualitative agreement. However, the displacements at which the fracture initiates is lower in the phase field approach. This can be accounted to the lower stresses in the region of fracture due to the degradation function~$g(d)$, cf.~\eqref{eq:Weak_form_deformation}, and to the higher observed values of $E_{frac}$ in both cases, before and after fracture propagation\footnote{From \eqref{eqn:FractureEnergy} it can be seen that $E_{frac}=0.5\,\text{N\,mm}$ before and $E_{frac}=1\,\text{N\,mm}$ after crack propagation is expected for a fully horizontal crack path with $G_c = 1\,\text{N/mm}$.}, which are presented in \Cref{fig:BVP2_phasefield_FractureEnergy} and \Cref{fig:BVP2_ShapeOpt_FractureEnergy} for both approaches.

%-----------------------------------------------------------------------

%=======================================================================
\section{Summary and conclusion}
\label{sec:SummaryConclusions}
%=======================================================================

We have outlined a commonly-used model for crack propagation and described how this model is adapted to a phase-field method and to a novel approach for the simulation of crack propagation based on shape optimization techniques. Here, the crack is interpreted as (possibly multiple) piecewise open curves which are then used as the argument of a minimization problem. We have employed gradient descent techniques on a shape space admitting an infinite dimensional \textsc{Riemann}ian manifold structure. This yields a mesh-independent approach and there is no requirement of staggered iterations. We have presented the numerical implementation of both, the phase-field and shape optimization approach, and applied them to two benchmark cases. We have analyzed the results and see good agreement for the fracture path in the first benchmark case with isotropic material behavior. Numerical results for second benchmark case are generated for anisotropic material behavior. More extensive analyses are performed for different anisotropy angles and also show identical directions of fracture between the two different approaches. Further, the boundary forces are compared and show some differences with respect to the displacement at which fracture initiates, which are then discussed.

%Future research includes the analysis of the effect of an angle-dependent failure criterion \TS{@CK/WW: should we write this, or is there too much risk of idea stealing?} and an extension of the approach to three-dimensional domains.

%=======================================================================
\section*{Acknowledgements}\label{Acknowledgements}
%=======================================================================
This research work is funded by the projects `Structural Health Monitoring' and
`KIBIDZ -- Intelligente Brandgefahrenanalyse für Geb\"{a}ude und Schutz der Rettungskr\"{a}fte durch K\"{u}nstliche Intelligenz und Digitale Brandgeb\"{a}udezwillinge'. Moreover,  `hpc.bw' also  provides computational resources (HPC cluster HSUper). All three projects are funded by dtec.bw -- Digitalization and Technology Research Center of the Bundeswehr.
dtec.bw is funded by the European Union -- NextGenerationEU.
%======================================================================================

%=======================================================================
\section*{Disclosure statement}

The authors report there are no competing interests to declare.

%=======================================================================

\bibliographystyle{plainurl} 
\bibliography{references.bib}


\end{document}